\documentclass{amsart}

\usepackage{a4wide,amssymb,tikz}
\usepackage{amssymb}
\usetikzlibrary{arrows}

\let\le\leqslant
\let\ge\geqslant
\let\x\times

\let\eps\varepsilon
\def\g{{\mathfrak g}}
\def\s{{\mathfrak s}}
\def\l{{\mathfrak l}}
\def\z{{\mathfrak z}}
\def\t{{\mathfrak t}}
\def\pr{^{\mathrm{pr}}}
\def\gl{{\mathfrak{gl}}}
\def\sl{{\mathfrak{sl}}}
\def\sp{{\mathfrak{sp}}}
\def\so{{\mathfrak{so}}}
\let\ro\rho

\def\abr#1{\langle#1\rangle}

\DeclareMathOperator{\diag}{diag}

\begin{document}

\title{On exceptional nilpotents in semisimple Lie algebras}
\author{A. G. Elashvili, V. G. Kac, E. B. Vinberg}

\maketitle

\section{Introduction}

\subsection{} In \cite{kw}, V. Kac and M. Wakimoto suggested
a construction of some class of rational vertex algebras coming from W-algebras. The data for this construction consists of a
positive integer $m$ and a nilpotent element $e$ of a semisimple Lie
algebra $\g$ satisfying some special conditions (see Section 2). A
pair $(m,e)$ satisfying those conditions was called exceptional by Kac
and Wakimoto. They classified exceptional pairs in simple Lie algebras
of type A$_n$. In this paper, we simplify the definition of
exceptional pairs and classify such pairs in all semisimple Lie
algebras. In particular, we prove that, for any semisimple Lie
algebra $\g$ and for any $m$, there is at most one, up to conjugation, nilpotent element $e$ in $\g$ such that the pair $(m,e)$ is exceptional.

\subsection{} Let $G$ be a connected semisimple algebraic group over
an algebraically closed field $F$ of characteristic 0, and let $\g=$Lie
$G$.

Let $\s\pr$ be a principal $\sl_2$-subalgebra of $\g$ \cite{kos}. The corresponding connected subgroup $S\pr\in G$ is isomorphic to $SL_2$ or $PSL_2$, in the former case its center being contained in the center of $G$. Let $h\pr$ be the semisimple element of $\s\pr$ corresponding to the matrix diag$(1,-1)\in \sl_2$.

For a positive integer $m$, let $\eps_{2m}\in F$ be a primitive $2m$-th root of 1, and let $s_m$ be the element of $S\pr$ corresponding to the matrix diag$(\eps_{2m},\eps_{2m}^{-1})$. Then $\sigma_m=Ad(s_m)$ is an automorphism of order $m$ of $\g$. We shall call $\sigma_m$ a principal automorphism of order $m$. Note that in general there are several conjugacy classes of principal automorphisms of order $m$ depending on the choice of $\eps_{2m}$.

The action of $\sigma_m$ can be described as follows. Choose a
maximal torus $T\subset G$ and a Weyl chamber in $\t=$Lie $T$
containing $h\pr$. Let $\{\alpha_1,...,\alpha_n\}$ be the corresponding set of simple roots. Then $\alpha_i(h\pr)=2$ for $i=1,\dots,n$. For any root $\alpha$, denote by $\mathrm{ht}(\alpha)$ its height, i.~e. the sum of coefficients in the linear expression of $\alpha$ in terms of $\alpha_1,...,\alpha_n$. Then, for any root vector $e_\alpha$,
$$
  \sigma_m(e_\alpha)=\eps_m^{\mathrm{ht}(\alpha)}e_\alpha,
$$
where $\eps_m=\eps_{2m}^2$.

For any automorphism $\sigma$ of $\g$, denote by $\g^\sigma$ the
subalgebra of fixed points of $\sigma$. Clearly, the subalgebra $\g^{\sigma_m}$ up to conjugacy does not depend on the choice of $\eps_{2m}$. Set
$$
  d(m)=\dim\g^{\sigma_m}.
$$

Recall that, for a simple Lie algebra $\g$, its Coxeter number
$h(\g)$ is defined as the order of the Coxeter element of the Weyl
group, and it is known \cite{kos} that
$$
  h(\g)=\mathrm{ht}(\delta)+1,
$$
where $\delta$ is the highest root of $\g$. The Coxeter numbers
of simple Lie algebras are given in the following table:

\

\begin{center}
\begin{tabular}{c|ccccccccc}
$\g$ & A$_n$ & B$_n$ & C$_n$ & D$_n$ & E$_6$ & E$_7$ & E$_8$ & F$_4$ & G$_2$\\
\hline
$h(\g)$ & $n+1$ & $2n$ & 2$n$ & $2n-2$ & 12 & 18 & 30 & 12 & 6
\end{tabular}
\end{center}

\

\noindent For any reductive Lie algebra $\l$, define the Coxeter number $h(\l)$ as the maximum of the Coxeter numbers of the simple factors of $\l$.
(For $\l$ abelian, set $h(\l)=1$.)

The above formula for the action of $\sigma_m$ implies that
$d(m)=$rk$(\g)$ (that is, $s_m$ is a regular element of $G$) if and
only if $m\ge h(\g)$. It is interesting that in the cases $m=h(\g)$ and $m=h(\g)+1$ the automorphism $\sigma_m$ is, up to conjugacy, the only regular inner automorphism of order $m$ \cite{kac2}.

\

\noindent{\bf Theorem 1}. (J.-P. Serre \cite{ser}.) \emph{Let $\sigma$
be an inner automorphism of $\g$ satisfying the condition
$\sigma^m=\mathrm{id}$ for a positive integer $m$. Then}
$$
 \dim \g^\sigma\ge d(m).
$$

\

A proof of this theorem has never been published. By a kind
permission of Serre, we include his elegant proof in Section 3 of
this paper.

\

\noindent{\bf Corollary}\cite{kac2}. \emph{Regular semisimple elements
$s\in G$ with $Ad(s)^m=\mathrm{id}$ exist if and only if $m\ge h(\g)$ (and in this case $s_m$ is one of such elements)}.

\

\subsection{} For a nilpotent element $e\in\g$, denote by $L(e)$ the centralizer of a maximal torus of the centralizer $Z(e)$ of $e$ in $G$. This is a (reductive) Levi subgroup of $G$ defined up to conjugacy by an
element of $Z(e)$. Its tangent Lie algebra $\l(e)$ contains $e$.
Define the Coxeter number $h(e)$ of $e$ as the Coxeter number of $\l(e)$.

{\noindent\bf Theorem 2}. \emph{For any nilpotent element $e\in\g$
with} $h(e)\le m$,
$$
  \dim Z(e)\ge d(m).
$$

\

This theorem looks similar to the above theorem of Serre. For even
nilpotent elements, it can be proved if making use of an unpublished
result of D.~Panyushev. But, unfortunately, we do not have a
conceptual proof of Theorem 2 in general case. In this paper, a proof
of it comes as a result of classification: see Sections 4,5.

\

{\noindent\bf Definition}. Let $m$ be a positive integer and $e\in\g$
be a nilpotent element. The pair $(m,e)$ is called
\emph{exceptional} if
\begin{align}
h(e)&\le m, \tag{I}\\
\dim Z(e)&=d(m). \tag{II}
\end{align}
\noindent In this case $e$ is called an exceptional nilpotent element
and $m$ an exceptional integer (for $\g$).
\

This definition does not coincide with but is essentially equivalent
to that of Kac and Wakimoto: see the discussion in Section 2.

Let $\g=\g_1+...+\g_s$ be the decomposition of $\g$ into a direct sum
of simple ideals, and let $e=e_1+...+e_s$ ($e_i\in\g_i$) be a
nilpotent element. Clearly, the pair $(m,e)$ is exceptional in $\g$
if and only if the pair $(m,e_i)$ is exceptional in $\g_i$ for every
$i$. Thus, the classification problem for exceptional pairs reduces
to the case, when $\g$ is simple.

The pair (1,0) is obviously exceptional. For all other exceptional
pairs, $m>1$ and $e\neq 0$. On the other hand, a pair $(m,e)$ with
$m\ge h(\g)$ is exceptional if and only if $e$ is regular (=principal).
The exceptional pairs of these two types are called \emph{trivial}.

A nilpotent element $e\in\g$ is said to be of \emph{principal type}
if it is principal in $\l(e)$. We shall prove (see Section 2) that any exceptional nilpotent element $e$ is of principal type.

\subsection{} The main result of this paper is the following classification
theorem.

\

\noindent{\bf Theorem 3}. \emph{For any simple Lie algebra $\g$ and
any positive integer $m$, there exists at most one nilpotent orbit $\mathrm{Ad}(G)e$ in $\g$ such that the pair $(m,e)$ is exceptional.
All non-trivial exceptional pairs $(m,e)$ in the classical and exceptional simple Lie algebras are listed in Tables 1 and 2, respectively.}

\

(In fact, Table 1 contains also some trivial exceptional pairs.)

In Table 1, the nilpotent element $e$ is given by the corresponding
partition of $N$, constituted by the orders of its Jordan blocks. In
Table 2, it is given by the type of the derived algebra $\l(e)'$ of
$\l(e)$. In the cases of G$_2$ and F$_4$, tildas mean that the root
system of the corresponding regular subalgebra consists of short roots.

Formulas for the numbers $d(m)$ in the classical case are given
in \ref{dmclass}. In the exceptional case, these numbers are given
in the tables on Figures~2-6.

\

\

\

\begin{center}
\ \hskip15em Table 1.\hskip17em Table 2.\hfill\

\

\

%\resizebox{!}{.35\textheight}
%{
\ \hfill \begin{tabular}{|c|c|l|}
\hline
  $\g$    &   $m$          &\hfill$e$\hfill\ \\
&&\\[-1em]
\hline
&&\\[-1em]
 $\sl_N$  &  any  &  $(m,...,m,r)$, $0\le r\le m-1$\\
&&\\[-1em]
\hline
 $\sp_N$, &  any  &  $(\underbrace{m,...,m}_{\textrm{even}},r)$, $0\le r\le m$, $r$
 even\\[-1em]
 $N$ even&&\\
&&\\[-1em]
\cline{2-3}
&&\\[-1em]
          &  odd  &
          $(\underbrace{m,...,m}_{\textrm{even}},m-1,m-1)$\\
&&\\[-1em]
\hline
&&\\[-1em]
 $\so_N$, & odd & $(\underbrace{m,...,m}_{\textrm{even}},r)$, $1\le r\le m$, $r$ odd\\[-1em]
 $N$ odd&&\\
&&\\[-1em]
\cline{2-3}
&&\\[-1em]
            & odd & $(\underbrace{m,...,m}_{\textrm{odd}},1,1)$\\
&&\\[-1em]
\cline{2-3}
&&\\[-1em]
            & even & $(m+1,\underbrace{m,...,m}_{\textrm{even}})$\\
&&\\[-1em]
\cline{2-3}
&&\\[-1em]
            & even&
            $(m+1,\underbrace{m,...,m}_{\textrm{even}},1,1)$\\
&&\\[-1em]
\cline{2-3}
&&\\[-1em]
            & even&  $(m+1,\underbrace{m,...,m}_{\textrm{even}},
            m-1,m-1)$\\
&&\\[-1em]
\hline
&&\\[-1em]
 $\so_N$,   & odd &  $(\underbrace{m,...,m}_{\textrm{even}},r,1)$, $1\le r\le m$, $r$ odd\\[-1em]
 $N$ even&&\\
&&\\[-1em]
\cline{2-3}
&&\\[-1em]
            & odd & $(\underbrace{m,...,m}_{\textrm{even}})$\\
&&\\[-1em]
\cline{2-3}
&&\\[-1em]
            &even &
            $(m+1,\underbrace{m,...,m}_{\textrm{even}},1)$\\
&&\\[-1em]
\cline{2-3}
&&\\[-1em]
            &even & $(m+1,\underbrace{m,...,m}_{\textrm{even}},m-1,m-1,1)$\\[1.5em]
\hline
\end{tabular}\hfill
%}
%
%\vfill
%
%\
%
%\pagebreak
%
%\
%
%Table 2.
%
%\
%
%\
%
\begin{tabular}{|c|c|l|}
\hline
  $\g$  &   $m$    &\hfill$e$\hfill\ \\
&&\\[-1em]
\hline
&&\\[-1em]
  G$_2$ &    2     & $\tilde{\mathrm A}_1$\\
&&\\[-1em]
\hline
&&\\[-1em]
  F$_4$ &    2     & A$_1+\tilde{\mathrm A}_1$\\
&&\\[-1em]
\cline{2-3}
&&\\[-1em]
        &    3     & $\tilde{\mathrm A}_2+$A$_1$\\
&&\\[-1em]
\hline
&&\\[-1em]
 E$_6$  &    2     & 3A$_1$\\
&&\\[-1em]
\cline{2-3}
&&\\[-1em]
        &    3     & 2A$_2$+A$_1$\\
&&\\[-1em]
\cline{2-3}
&&\\[-1em]
        &    5     & A$_4$+A$_1$\\
&&\\[-1em]
\cline{2-3}
&&\\[-1em]
        &    8     & D$_5$\\
&&\\[-1em]
\hline
&&\\[-1em]
 E$_7$  &   2     & 4A$_1$\\
&&\\[-1em]
\cline{2-3}
&&\\[-1em]
        &   3     & 2A$_2$+A$_1$\\
&&\\[-1em]
\cline{2-3}
&&\\[-1em]
        &   4     & A$_3$+A$_2$+A$_1$\\
&&\\[-1em]
\cline{2-3}
&&\\[-1em]
        &  5    &  A$_4$+A$_2$\\
&&\\[-1em]
\cline{2-3}
&&\\[-1em]
        &  7    &  A$_6$\\
&&\\[-1em]
\hline
&&\\[-1em]
  E$_8$ &    2  &   4A$_1$\\
&&\\[-1em]
\cline{2-3}
&&\\[-1em]
        & 3     & 2A$_2$+2A$_1$\\
&&\\[-1em]
\cline{2-3}
&&\\[-1em]
        & 4     & 2A$_3$\\
&&\\[-1em]
\cline{2-3}
&&\\[-1em]
        & 5     & A$_4$+A$_3$\\
&&\\[-1em]
\cline{2-3}
&&\\[-1em]
        & 7     & A$_6$+A$_1$\\
&&\\[-1em]
\cline{2-3}
&&\\[-1em]
        & 8     & A$_7$\\
\hline
\multicolumn{3}{c}{}\\[-.2em]
\multicolumn{3}{c}{}\\
\multicolumn{3}{c}{}\\
\multicolumn{3}{c}{}
\end{tabular}

\end{center}

\

\

\subsection{}\label{un} If $(m,e)$ is an exceptional pair, then $m\ge h(e)$ and
$$
  d(h(e))\le\dim Z(e)=d(m)\le d(h(e)),
$$
whence $d(m)=d(h(e))$ and $(h(e),e)$ is also an exceptional pair. A priori it is possible that $m>h(e)$. This really happens for regular $e$, where any $m>h(e)=h(\g)$ fits. Apart from this case, this never happens in the exceptional Lie algebras. In the classical Lie algebras, this happens only in the following two cases:
\begin{itemize}
\item[1)] $\g=\sp_{2n}$, $n$ even, $e$ is defined by the partition $(n,n)$, $m=n+1$ ($>n=h(e)$);
\item[2)] $\g=\so_{2n+1}$, $e$ is defined by the partition $(2n-1,1,1)$, $m=2n-1$ ($>2n-2=h(e)$).
\end{itemize}

\subsection{} It follows from the tables in \cite{a} that the centralizer of any exceptional nilpotent element in a simple Lie algebra $\g\ne\so_N$ is connected. For $\g=\so_N$, it may have at most two connected components.

\subsection{} For any positive integer $k$ and any simple Lie algebra $\g$, set $N_k(\g)=\{x\in\g : (\mathrm{ad}x)^k=0\}$. The irreducible components of the varieties $N_k(\g)$ are found in \cite{v}. One can check that any exceptional nilpotent orbit is open in some $N_k(\g)$.

\subsection{} We thank D.~Panyushev who let us know about an unpublished result of J.-P.~Serre (Theorem 1 above) having been mentioned in Serre's talk on an Oberwolfach conference in 1998, and to J.-P.~Serre who permitted us to include his proof in this paper and made some useful remarks on the paper. We also thank M.~Jibladze for providing technical help in the preparation of this manuscript.

This work was mainly done during the stay of the first and third authors
at the University of Bielefeld in July of 2008, supported by SFB 701. We thank this university for its hospitality. The first author also acknowledges partial support from GNSF (Grant \# ST07/3-174) and SFB/TR 12 of the Deutsche Forschungsgemeinschaft.

\section{Definition of exceptional pairs}

\subsection{} The original definition of exceptional pairs was given in terms
of W-algebras. However, Theorem 2.32 of \cite{kw} permits to give an
equivalent definition in internal terms of the algebra $\g$.

For a nilpotent element $e$ of a semisimple Lie algebra $\g$ and a
positive integer $m$, denote by $S(m,e)$ the set of all regular
semisimple elements $s$ of $L(e)$ such that  Ad$(s)^m$=id. According
to \cite{kw}, the pair $(m,e)$ is exceptional if $e$ is of
principal type, $S(m,e)\ne\varnothing$, and
$$
  \min_{s\in S(m,e)}\dim Z(s)=\dim Z(e).
$$

\

\noindent{\bf Remark}. In fact, Kac and Wakimoto require in addition
that $m$ should be coprime to the "lacety" of $\g$, which is 1 for
types A,D,E, 2 for B,C,F, and 3 for G. We will disregard this
requirement, which is natural from the point of view of $W$-algebras but
looks artificial from the point of view of the theory of semisimple Lie algebras and, besides, it does not facilitate the classification. One can note, however, that this requirement is violated in cases 1) and 2) of
subsection \ref{un}, so if we adopt it, then for any non-principal nilpotent element $e\in\g$ there will be at most one positive integer $m$ such that the pair $(m,e)$ is exceptional, as was conjectured in \cite{kw}.

\

Applying Corollary to Theorem 1 to $L(e)$, we obtain that
$S(m,e)\ne\varnothing$ if and only if $h(e)\le m$, which is just
condition (I) of the definition of an exceptional pair given in the
introduction.

Let us fix a maximal torus $T\subset G$ and a set of simple roots
$\alpha_1$,..., $\alpha_n$ with respect to it. Fix also a principal
$\sl_2$-subalgebra $\s\pr$ with $h\pr$ contained in the Weyl chamber
in $\t=$Lie $T$.
Let $s_m$ be the element of the corresponding subgroup $S\pr$ defined as in the introduction. One may assume that $L(e)$ contains $T$ and, moreover, that $\l(e)$ is generated by $\t=$Lie $T$ and some positive and the opposite negative simple root vectors. Then $s_m\in T\subset L(e)$, and the description of the action of $\sigma_m=$Ad($s_m$) given in the introduction shows that if $h(e)\le m$, then $s_m$ is a regular element of $L(e)$, so $s_m\in S(m,e)$. Now, Theorem 1 implies that
$$
  \min_{s\in S(m,e)}\dim Z(s)=\dim Z(s_m)=d(m).
$$
Thus, the dimension condition in the above definition of an
exceptional pair reduces to condition (II) of the definition given
in the introduction.

\subsection{} Let $e\in\g$ be a nilpotent element with $h(e)\le m$, and let
$e_0$ be a principal nilpotent element of $\l(e)$. Then
$\l(e_0)=\l(e)$, so $h(e_0)=h(e)\le m$. Further, $e$ lies in the
closure of the $L(e)$-orbit of $e_0$ and, the more, in the closure
of the $G$-orbit of $e_0$. Hence,
$$
  \dim Z(e)\ge \dim Z(e_0),
$$
the equality taking place only if $e\in $Ad$(G)e_0$.

According to Theorem 2 (which will be proved in Sections 4,5 together
with the classification of exceptional pairs), dim $Z(e_0)\ge d(m)$.
Hence, the equality dim $Z(e)=d(m)$ can only take place if $e\in
$Ad$(G)e_0$.

Suppose that $e=$Ad$(g)e_0$ for some $g\in G$. Multiplying $g$ from
the left by some element of $Z(e)$, one may assume that
Ad$(g)\l(e)=\l(e)$. Then Ad($g$) leaves invariant the principal
nilpotent orbit in $\l(e)$. Hence, $e$ lies in this orbit, i.e., $e$
is a principal nilpotent element in $\l(e)$.

Thus, the original definition of an exceptional pair given in
\cite{kw} is in fact equivalent to the definition given in the
introduction.

\section{Proof of Theorem 1}

\subsection{} Serre's proof of Theorem 1 is based on the product
formula for the character $\varphi_\lambda$ of the restriction
to $S\pr$ of the irreducible representation of $G$ with highest
weight $\lambda$. Apparently, this formula was already known in
the 60s but was not explicitly written at that time. A more recent
reference is \cite{kac1}, formula (3.29).

Denote by $s(t)$ the element of $S\pr$ corresponding to the matrix
diag$(t,t^{-1})\in$SL$_2$. Then the formula is
$$
\varphi_\lambda(s(t))=t^{-\abr{\lambda,
\ro^\vee}}\prod_{\alpha>0}\frac{t^{\abr{\lambda+\ro,\alpha^\vee}}-1}{t^{\abr{\ro,\alpha^\vee}}-1},
$$
where $\ro$ is, as usual, the half-sum of positive roots,
$\alpha^\vee$ denotes the coroot corresponding to $\alpha$, and
$\ro^\vee$ is the half-sum of positive coroots. It is obtained from
the Weyl character formula; one should only note that the Weyl
denominator formula (for the dual root system) is applicable to the
numerator if one is only interested in the restriction of the
character to $S\pr$.

Clearly, the product in the right hand side must be a polynomial in
$t$. Let us represent it in a little different form. First, note
that $\abr{\ro,\alpha^\vee}$ is nothing else than the height
ht$(\alpha^\vee)$ of $\alpha^\vee$, since $\abr{\ro,\alpha_i^\vee}=1$
for every simple coroot $\alpha_i^\vee$. Second,
$\abr{\lambda+\ro,\alpha^\vee}$ can be interpreted as the ``weighted
height'' of $\alpha$ if one assigns the weight
$l_i=\abr{\lambda+\ro,\alpha_i^\vee}=\abr{\lambda,\alpha_i^\vee}+1$ to
each simple coroot $\alpha_i^\vee$. Denote the so defined weighted
height of $\alpha^\vee$ by ht$_{\boldsymbol l}(\alpha^\vee)$, where
$\boldsymbol l=(l_1,...,l_n)$. Note that $l_1,...,l_n$ may be arbitrary positive integers. In particular, ht$(\alpha^\vee)=$ht$_{\boldsymbol1}(\alpha^\vee)$, where
$\boldsymbol1=(1,...,1)$.

Finally, replacing the root system of $G$ with its dual, we
come to the following

\

\noindent{\bf Proposition 2.1} \cite[Theorem 1]{sta}. \emph{For any set
$\boldsymbol l=(l_1,...,l_n)$ of positive integers, the polynomial
$$
  P_{\boldsymbol l}(t)=\prod_{\alpha>0}{(t^{\mathrm{ht}_{\boldsymbol l}(\alpha)}-1)}
$$
is divisible by the polynomial}
$$
  P_{\boldsymbol1}(t)=\prod_{\alpha>0}{(t^{\mathrm{ht}(\alpha)}-1)}.
$$

\subsection{} Now we are ready to prove Theorem 1.

Let $\sigma$ be an inner automorphism of $\g$ satisfying the
condition $\sigma^m=$id for some positive integer $m$. One may
assume that $\sigma$ is a conjugation by some element $s\in T$.
Then
$$
  \sigma(e_{\alpha_i})=\eps_m^{l_i}e_{\alpha_i}\,(l_i\in\{1,...,m\}),
$$
where $\eps_m=\exp\frac{2\pi i}m$, and, hence, for any $\alpha>0$,
$$
  \sigma(e_\alpha)=\eps_m^{\mathrm{ht}_{\boldsymbol l}(\alpha)}e_\alpha,
$$
where $\boldsymbol l=(l_1,...,l_n)$. This implies that
$$
  \dim\g^\sigma=\mathrm{rk}\,\g+2\#\{\alpha>0: m\,|\,\mathrm{ht}_{\boldsymbol l}(\alpha)\}.
$$
In particular,
$$
  d(m)=\mathrm{rk}\,\g+2\#\{\alpha>0: m\,|\,\mathrm{ht}(\alpha)\}.
$$

Clearly, $\#\{\alpha>0: m\,|\,\mathrm{ht}(\alpha)\}$ is the multiplicity
of $\eps_m$ as a root of the polynomial $P_{\boldsymbol1}(t)$, while
$\#\{\alpha>0: m\,|\,\mathrm{ht}_{\boldsymbol l}(\alpha)\}$ is the multiplicity of
$\eps_m$ as a root of the polynomial $P_{\boldsymbol l}(t)$. According to
Proposition 2.1, the latter is not less than the former, whence
Theorem 1 follows.

\subsection{} In addition, let us prove some useful monotonicity properties
of the function $m\mapsto d(m)$.

\

\noindent{\bf Proposition 2.2}. \emph{If $m'<m$, then $d(m')\ge
d(m)$. Moreover, if $m'\,|\,m$, $m'\neq m<h(\g)$, then $d(m')>d(m)$.}

\begin{proof} The second inequality immediately follows from the
preceding formula for $d(m)$ and the definition of $h(\g)$.

The first inequality follows from the well-known fact that the
number $\#\{\alpha>0:\mathrm{ht}(\alpha)=k\}$ is monotonically
decreasing in $k$ (see \cite{dyn,kos}). Indeed, if $m'<m$, then, for
any $k$,
$$
  \#\{\alpha>0:\mathrm{ht}(\alpha)=km'\}\ge\#\{\alpha>0:\mathrm{ht}(\alpha)=km\},
$$
whence the required inequality follows.
\end{proof}

\section{Classification: the classical Lie algebras}

\subsection{} Our strategy in proving Theorems 2 and 3 will be the following. For each simple Lie algebra $\g$ we consider the set Nil$(\g)$ of its nilpotent orbits partially ordered by the inclusion of the closures. Then for each $m<h(\g)$ we consider the subset Nil$_m(\g)$ of nilpotent orbits Ad$(G)e$ with $h(e)\le m$ and determine its maximal elements, which we call \emph{essential} nilpotent orbits. (They are automatically nilpotent orbits of principal type: see Section 2.) We check that dim$Z(e)\ge d(m)$ for any essential nilpotent orbit and thereby prove Theorem 2. At the same time, we find all essential nilpotent orbits Ad$(G)e\in$Nil$_m(\g)$ with dim $Z(e)=d(m)$ and thus obtain a classification of exceptional pairs. It turns out that for each $m$
there is at most one essential nilpotent orbit with this property.

\subsection{} First of all, we will deduce some general formulas
for the dimensions of the centralizers of semisimple elements in the classical groups $G=SL_N$, $Sp_N$, $SO_N$.

Let $s\in SL_N$ be a semisimple element with eigenvalues of
multiplicities $n_1,...,n_p$ (so $n_1+...+n_p=N$). Denote by $K$ the sum
of squares of these multiplicities. The centralizer of $s$ in $GL_n$ is
isomorphic to $GL_{n_1}\x...\x GL_{n_p}$ and, hence, its dimension is
equal to $K$. The dimension of the centralizer $Z(s)$ of s in $SL_N$ is
one less. Thus,
$$
  \dim Z(s)=K-1\textrm{  for }G=SL_N.
\eqno{(1)}
$$

If $s\in Sp_N$ or $SO_N$, the eigenvalues of $s$ distinct from $\pm1$
decompose into pairs of mutually inverse ones. Let $n_1,...,n_q$ be
the common multiplicities of the eigenvalues of these pairs, and let
$n_+$ and $n_-$ be the multiplicities of the eigenvalues $1$ and $-1$ (so
$2(n_1+...+n_q)+n_++n_-=N)$. As above, denote by $K$ the sum of squares
of all the multiplicities, that is,
$$
  K=2(n_1^2+...+n_q^2)+n_+^2+n_-^2.
$$

For $s\in Sp_N$ the centralizer of $s$ in $Sp_N$ is isomorphic to $GL_{n_1}\x...\x GL_{n_q}\x Sp_{n_+}\x Sp_{n_-}$. Hence,
$$
  \dim Z(s)=n_1^2+...+n_q^2+\frac{n_+(n_++1)}{2}+\frac{n_-(n_-+1)}{2},
$$
which can be written in the form
$$
  2\dim Z(s)=K+n_++n_-\textrm{  for }G=Sp_N.
\eqno{(2)}
$$

Similarly, for $s\in SO_N$ the connected centralizer of $s$ in $SO_N$ is
isomorphic to $GL_{n_1}\x...\x GL_{n_q}\x SO_{n_+}\x SO_{n_-}$, whence
$$
2\dim Z(s)=K-n_+-n_-\textrm{  for }G=SO_N.
\eqno{(3)}
$$

\subsection{}\label{dmclass} Let us now calculate the numbers $d(m)$ for the classical simple Lie algebras
$\g=\sl_N$, $\sp_N$, $\so_N$. In the last two cases, we will suppose that
the invariant (skew-symmetric or symmetric) inner product is defined by
\begin{align*}
  (e_i,e_{N+1-i})&=1\textrm{ for }i\le (N+1)/2,\\
  (e_i,e_j)&=0\textrm{ for }i+j\neq (N+1),
\end{align*}
where $\{e_1,\dots, e_N\}$ is the standard basis of $F^N$.

Fix a maximal torus $T$ in $G$ consisting of diagonal matrices and a Borel subgroup consisting of upper triangular matrices. If $\g\neq \so_N$ with $N$ even, then
$$
  h\pr=\diag(N-1,N-3,...,-(N-3),-(N-1)),
$$
and, hence,
$$
  s_m=\diag(\eps_{2m}^{N-1},\eps_{2m}^{N-3},...,\eps_{2m}^{-(N-3)},\eps_{2m}^{-(N-1)}).
\eqno{(4)}
$$
If $\g=\so_N$ with $N$ even, then the subgroup $S\pr$ is contained in the subgroup $SO_{N-1}$ embedded into $SO_N$ in the standard way (and is a principal 3-dimensional subgroup there). It follows that in this case
$$
  s_m=\diag(\eps_{2m}^{N-2},\eps_{2m}^{N-4},...,\eps_{2m}^2,1,1,\eps_{2m}^{-2},...,\eps_{2m}^{-(N-4)},\eps_{2m}^{-(N-2)}).
\eqno{(5)}
$$

Let $N=qm+r$, where $1\le r\le m$, if $\g=\so_N$, $N$ even, and $0\le r\le m-1$ in
all the other cases. Then $(\underbrace{m,...,m}_q,r)$ is a partition of $N$. Denote by $K(m)$ the sum of squares of the parts of the dual partition $(\underbrace{q+1,...,q+1}_r,\underbrace{q,...,q}_{m-r})$, that is,
$$
  K(m)=r(q+1)^2+(m-r)q^2.
$$

\

\noindent{\bf Proposition 3.1.}
\emph{
\begin{itemize}
\item[1)] For $\g=\sl_N$,
$$
  d(m)=K(m)-1.
$$
\item[2)] For $\g=\sp_N$,
$$
2d(m)=K(m)+ \begin{cases}
              q,&\textrm{if $m$ is odd, $q$ is even,}\\
  q+1,&\textrm{if $m$ and $q$ are odd,}\\
              0,&\textrm{if $m$ is even.}
\end{cases}
$$
\item[3)] For $\g=\so_N$, $N$ odd,
$$
2d(m)=K(m)- \begin{cases}
              q,&\textrm{if $m$ and $q$ are odd,}\\
  q+1,&\textrm{if $m$ is odd, $q$ is even,}\\
              2q+1,&\textrm{if $m$ is even.}
\end{cases}
$$
\item[4)] For $\g=\so_N$, $N$ even,
$$
  2d(m)=K(m)- \begin{cases}
q,&\textrm{if $m$ is odd, $q$ is even,}\\
              q+1,&\textrm{if $m$ and $q$ are odd,}\\
2q, &\textrm{if $m$ and $q$ are even,}\\
              2(q+1),&\textrm{if $m$ is even, $q$ is odd.}\\
\end{cases}
$$
\end{itemize}
}

\begin{proof}
The proof of 1)-3) is obtained by applying formulas (1)-(3) to $s=s_m$.
Since the eigenvalues of $s_m$ constitute a geometric progression with denominator $\eps_m$, their multiplicities are $\underbrace{q+1,...q+1}_r,\underbrace{q,...,q}_{m-r}$. In particular, $K=K(m)$, whence 1) immediately follows.

To prove 2) and 3), one should determine $n_++n_-$ for $s=s_m$.

For $\g=\sp_N$, it follows from (4) that all the eigenvalues of $s_m$ are $m$-th
roots of -1. Hence, 1 is not an eigenvalue of $s_m$. Moreover, if $m$ is even,
-1 is not an eigenvalue, neither. If $m$ is odd, $n_+=q$ or $q+1$. In order to distinguish between these two possibilities, it suffices to note that for symmetry reason $n_+$ must be even.

For $\g=\so_N$, $N$ odd, the eigenvalues of $s_m$ are $m$-th roots of 1. Hence, if $m$
is odd, -1 is not an eigenvalue, while $n_+=q$ or $q+1$; but for symmetry reason $n_+$
must be odd, which permits to determine $n_+$ uniquely. If $m$ is even, both $n_+$ and $n_-$ are equal to $q$ or $q+1$. For symmetry reason, $n_++n_-$ must be odd, whence $n_++n_-=2q+1$.

Let now $\g=\so_N$ with $N$ even. It follows from our definition of $K(m)$ that $K=K(m)$ if $q$ is odd, and $K=K(m)+2$ if $q$ is even. As in the preceding case, the eigenvalues of $s_m$ are $m$-th roots of 1 (see (5)). If $m$ is odd, -1 is not an eigenvalue, while the multiplicity of the eigenvalue 1 is even and equals $q+1$ or $q+2$. If $m$ is even, we have $n_+=q+1$ or $q+2$ and $n_-=q$ or $q+1$; but for symmetry reason $n_+$ and $n_-$ are even, so $n_++n_-=2q+2$, which gives 4).
\end{proof}

\subsection{} In this subsection, we collect some well-known facts about nilpotent orbits in the classical simple Lie algebras $\g=\sl_N$, $\sp_N$, $\so_N$. For more details and proofs, see, for example, \cite{cm}.

A nilpotent orbit Ad$(G)e$ in $\g$ is uniquely defined by the partition $(n_1,...n_p)$ of $N$ constituted by the orders of Jordan blocks of $e$
(acting on $F^N$), with the only reservation that in the case $\g=\so_N$ with $N\equiv0$~(mod 4), the partitions with all even parts correspond to two different nilpotent orbits permuted by an outer automorphism of $\g$. The partition $(n_1,...,n_p)$ may be arbitrary for $\g=\sl_N$ but in the other cases is subject to some restrictions. Namely, for $\g=\sp_N$ the multiplicity of each odd part of the partition should be even, while for $\g=\so_N$ the multiplicity of each even part should be even; we shall call such partitions \emph{admissible} for $\g$. We agree to think of the parts of a partition as of the rows of a Young diagram going from the bottom to the top and aligned from the left.

Denote by $K(e)$ the sum of squares of the parts of the partition dual to $(n_1,...,n_p)$ (constituted by the columns of the corresponding Young diagram). Then the dimension of the centralizer $Z(e)$ of $e$ in $G$ is given by the following formulas:
$$
  \dim Z(e)=K(e)-1\textrm{  for  }\g=\sl_N,
\eqno{(6)}
$$
$$
  2\dim Z(e)=K(e)+\#\{i: n_i \textrm{ odd}\}\textrm{  for  }\g=\sp_N,
\eqno{(7)}
$$
$$
  2\dim Z(e)=K(e)-\#\{i: n_i \textrm{ odd}\}\textrm{  for }\g=\so_N.
\eqno{(8)}
$$

To describe the partial order on the set Nil($\g$) of nilpotent orbits in $\g$, let us introduce the notion of a ``\emph{simple crumbling}'' of a partition $(n_1,...,n_p)$ as
the transition to a partition of the form
$$
  (n_1,...n_{i-1},n_i+1,n_{i+1},...,n_{j-1},n_j-1,n_{j+1},...,n_p),
$$
provided $n_{i-1}>n_i$ and $n_j>n_{j+1}$. For example, on Fig. 1 the simple crumbling of the partition (8,6,6,3,2) to the partition (8,7,6,2,2) is shown.

\

\

\begin{center}
\begin{tikzpicture}[scale=.5]
\fill[gray!20] (2,3) rectangle (3,4);
\draw (1,0) -- (1,5);
\draw (2,0) -- (2,4);
\draw (3,0) -- (3,3);
\draw (4,0) -- (4,3);
\draw (5,0) -- (5,3);
\draw (6,0) -- (6,1);
\draw (7,0) -- (7,1);
\draw (0,1) -- (6,1);
\draw (0,2) -- (6,2);
\draw (0,3) -- (3,3);
\draw (0,4) -- (2,4);
\draw[dashed] (7,1) -- (7,2) -- (6,2);
\draw[thick] (0,0) -- (8,0) -- (8,1) -- (6,1) -- (6,3) -- (3,3) -- (3,4) -- (2,4) -- (2,5) -- (0,5) -- cycle;
\draw[-stealth] (2.5,3.5) .. controls (6.5,3.5) and (6.5,3.5) .. (6.5,1.5);
\node at (11,2.5) {$\longmapsto$};
\fill[gray!20] (20,1) rectangle (21,2);
\draw (15,0) -- (15,5);
\draw (16,0) -- (16,3);
\draw (17,0) -- (17,3);
\draw (18,0) -- (18,3);
\draw (19,0) -- (19,3);
\draw (20,0) -- (20,2);
\draw (21,0) -- (21,1);
\draw (14,1) -- (21,1);
\draw (14,2) -- (20,2);
\draw (14,3) -- (16,3);
\draw (14,4) -- (16,4);
\draw[dashed] (17,3) -- (17,4) -- (16,4);
\draw[thick] (14,0) -- (22,0) -- (22,1) -- (21,1) -- (21,2) -- (20,2) -- (20,3) -- (16,3) -- (16,5) -- (14,5) -- cycle;
\end{tikzpicture}

\

Fig.1

\end{center}

\

\

Let Ad$(G)e$ and Ad$(G)e'$ be two nilpotent orbits in $\g$ corresponding to partitions
$(n_1,...,n_p)$ and $(n'_1,...,n'_{p'})$. Then Ad$(G)e$ lies in the closure of Ad$(G)e'$ if and only if the partition $(n'_1,...,n'_{p'})$ can be obtained from $(n_1,...n_p)$ by consecutive simple crumblings (without assuming that all the intermediate partitions should be admissible) (see \cite{cm}).

\subsection{} As was explained in Section 2, for our purposes it suffices to consider only nilpotent elements of principal type. Let us describe such elements in terms of partitions (cf. \cite{kw}).

First of all, in all the classical simple Lie algebras, but $\so_N$ with $N$ even, a principal nilpotent element is defined by the trivial partition $(N)$. In $\so_N$ with N even, it is defined by the partition $(N-1,1)$.

Let $e\in\sl_N$ be a nilpotent element defined by a partition $(n_1,...,n_p)$. Then e is conjugate to a principal nilpotent element of the Levi subalgebra consisting of the matrices
$$
  A=\diag(A_1,...,A_p)\  (A_1\in\gl_{n_1},...,A_p\in\gl_{n_p})
$$
with tr $A=0$. Thus, all nilpotent elements in $\gl_N$ are of principal type.

In $\g=\sp_N$ or $\so_N$, any Levi subalgebra $\l$ consists of the matrices of the form
\begin{multline*}
  A=\diag(A_1,...,A_s,A_0,-A'_s,...,-A'_1)\\  (A_1\in\gl_{n_1},..., A_s\in\gl_{n_s}, \       A_0\in\sp_{n_0}\textrm{ or }\so_{n_0}, \textrm{ resp.; } (2(n_1+...+n_s)+n_0=N),
\end{multline*}
where $'$ denotes the transposition with respect to the second diagonal. Thus, $\l$ is isomorphic to $\gl_{n_1}+...+\gl_{n_s}+\sp_{n_0}$ or $\gl_{n_1}+...+\gl_{n_s}+\so_{n_0}$, resp. This implies the following characterization of nilpotent elements of principal type in terms of the corresponding partition $(n_1,...,n_p)$:

\begin{itemize}
\item[1)] for $\g=\sp_N$, the multiplicities of all parts of the partition, except for at most one even part, should be even;
\item[2)] for $\g=\so_N$, $N$ odd, the multiplicities of all parts of the partition, except for at most one odd part, should be even;
\item[3)] for $\so_N$, $N$ even, either all the multiplicities are even, or the multiplicities of 1 and some other odd part are odd, while all the other multiplicities are even.
\end{itemize}

We shall refer to such partitions as to (admissible) \emph{partitions of principal type}.

It follows from this description and the table of the Coxeter numbers of simple Lie algebras (see the introduction) that, for $\g =\sl_N$ or $\sp_N$, the Coxeter number of a nilpotent element of principal type corresponding to the partition $(n_1,...,n_p)$, is equal to $n_1$ (the maximum of the parts of the partition). For $\g=\so_N$, it is equal to $n_1$ or $n_1-1$, the latter taking place iff $n_1$ is odd and $n_1>n_2$.

\subsection{} Case $\g=\sl_N$. In this case, Nil$_m(\g)$ consists of the nilpotent orbits defined by the partitions all whose parts do not exceed $m$. Any such partition crumbles to the partition $(m,...,m,r)$ with $0\le r\le m-1$, which is thereby the only maximal element of Nil$_m(\g)$. Let Ad$(G)e$ be the corresponding nilpotent orbit. Then, by (6) and Proposition 3.1.1),
$$
  \dim Z(e)=r(q+1)^2+(m-r)q^2-1=d(m).
$$
This proves Theorems 2 and 3 for $\sl_N$ (cf. \cite{kw}).

\subsection{} Case $\g=\sp_N$. In this case, the nilpotent orbits of principal type in Nil$_m(\g)$ are defined by the partitions of principal type all whose parts do not exceed $m$. Any such partition crumbles to one of the following partitions of the same class:
\begin{itemize}
\item[1)] $(\underbrace{m,...,m}_{\textrm{even}},r)$ with $0\le r\le m$;
\item[2)] $(\underbrace{m,...,m}_{\textrm{even}},s,s)$ with $\frac m2+1\le s\le m-1$;
\item[3)] $(\underbrace{m,...,m}_{\textrm{even}},m-1,s,s)$ with $m$ odd, $1\le s\le\frac{m-3}2$;
\item[4)] $(\underbrace{m,...,m}_{\textrm{odd}},s,s)$ with $m$ even, $1\le s\le m-1$.
\end{itemize}

Applying (7) and Proposition 3.1.2) to partitions 1)-4), we obtain Theorems 2 and 3 for $\sp_N$. The calculations can be simplified if one notes that the difference $\dim Z(e)-d(m)$ does not change when deleting $2k$ parts equal to $m$ from the partition (and diminishing $N$ by $2km$). In case 1) this reduces the consideration to the partition $(r)$, which corresponds to the principal nilpotent orbit in $\sp_r$; hence, in this case the pair $(m,e)$ is always exceptional. In case 2) it suffices to consider the partition $(s,s)$, where we obtain
$$
\begin{array}{rl}
2(\dim Z(e)-d(m))
&=\Bigg[4s+\begin{cases} 0, &\textrm{if $s$ is even},\\2, &\textrm{if $s$ is odd}\end{cases}\Bigg]
-\Bigg[(6s-2m)+\begin{cases} 0, &\textrm{if $m$ is even},\\2, &\textrm{if $m$ is odd}\end{cases}\Bigg]\\
&=2(m-s)+\begin{cases} 0, &\textrm{if $s$ is even},\\2, &\textrm{if $s$ is odd}\end{cases}
-\begin{cases} 0, &\textrm{if $m$ is even},\\2, &\textrm{if $m$ is odd}\end{cases}
\ \ \ \ge 0,
\end{array}
$$
the equality taking place iff $m$ is odd and $s=m-1$. The cases 3) and 4) are treated similarly.

\subsection{} Case $\g=\so_N$, $N$ odd. In this case, the nilpotent orbits of principal type in Nil$_m(\g)$ are defined by the partitions of principal type all whose parts do not exceed $m$ or, if $m$ is even, by the partitions of principal type, whose maximal part is
equal to $m+1$ and occurs with multiplicity 1. Any such partition crumbles to one of the following partitions of the same class:
\begin{itemize}
\item[1)] $(\underbrace{m,...,m}_{\textrm{even}},r)$ with $1\le r\le m$;
\item[2)] $(\underbrace{m,...,m}_{\textrm{odd}},s,s)$ with $m$ odd, $1\le s\le\frac{m-1}2$;
\item[3)] $(\underbrace{m,...,m}_{\textrm{even}},s,s,1)$ with $m$ odd, $\frac{m+1}2\le s\le m-1$;
\item[4)] $(m+1,\underbrace{m,...,m}_{\textrm{even}},s,s)$ with $m$ even, $0\le s\le m-1$.
\end{itemize}

Applying (8) and Proposition 3.1.3) to partitions 1)-4), we obtain Theorems 2 and 3  for $\so_N$, $N$ odd. To simplify the calculations, one can note that in cases 1)-3) the difference dim$Z(e)-d(m)$ does not change when deleting $2k$ parts equal to $m$ from the partition, if $m$ is odd, and decreases, if $m$ is even. In case 1) this reduces the consideration to the partition $(r)$, which corresponds to the principal nilpotent orbit in $\so_r$; hence, in this case the pair $(m,e)$ is always exceptional if $m$ is odd, while if $m$ is even, it is exceptional only in the trivial case when $N=r$ (and, hence, $m>h(e)$).

In case 4), a direct calculation shows that for $s<\frac m2$
$$
  2(\dim Z(e)-d(m))=2s-\begin{cases} 0,&\textrm{if $s$ is even},\\
                       2,&\textrm{if $s$ is odd}
\end{cases}\ \ \ \ge 0,
$$
the equality taking place iff $s=0$ or 1. For $s\ge\frac m2$
$$
  2(\dim Z(e)-d(m))=2(m-s)-\begin{cases} 0,&\textrm{if $s$ is even},\\
                           2,&\textrm{if $s$ is odd}
\end{cases}\ \ \ \ge 0,
$$
the equality taking place iff $s=m-1$.

\subsection{} Case $\g=\so_N$, $N$ even. The partitions defining nilpotent orbits of principal type in Nil$_m(\g)$ are described in the same way as for $\so_N$, $N$ odd. Any such partition crumbles to one of the following partitions of the same class:
\begin{itemize}
\item[1)] $(\underbrace{m,...,m}_{\textrm{even}},r,1)$ with $1\le r\le m$;
\item[2)] $(\underbrace{m,...,m}_{\textrm{even}},s,s)$ with $\frac m2<s\le m$;
\item[3)] $(\underbrace{m,...,m}_{\textrm{odd}},s,s,1)$ with $m$ odd, $1\le s\le\frac{m-1}2$;
\item[4)] $(m+1,\underbrace{m,...,m}_{\textrm{even}},s,s,1)$ with $m$ even, $1\le s\le m$.
\end{itemize}
Applying (8) and Proposition 3.1.4) to partitions 1)-4), we obtain Theorems 2 and 3 for $\so_N$, $N$ even. The calculations are similar to those in the preceding case.

\section{Classification: the exceptional Lie algebras}

\subsection{} To classify exceptional pairs in the exceptional simple Lie algebras, we need a  formula for computing the numbers $d(m)$.
Let $\g$ be a simple Lie algebra of rank $n$, and $\s\pr$ be a principal $\sl_2$-subalgebra of $\g$.
It is well-known \cite{kos}
that the adjoint representation of $\s\pr$ in $\g$
decomposes into a sum of $n$ irreducible representations
of dimensions $2m_1+1,...,2m_n+1$,
where $m_1,...,m_n$ are the exponents of $\g$.
Clearly, the eigenspace of ad($h\pr$) corresponding to the eigenvalue $2k>0$ is spanned by the positive root vectors $e_\alpha$ with ht$(\alpha)=k$.
It follows that
$$
  \#\{\alpha>0: \mathrm{ht}(\alpha)=k\}=\#\{i: m_i\ge k\},
$$
Hence,
$$
  \#\{\alpha>0: m|\mathrm{ht}(\alpha)\}=\sum_{i=1}^n\left[\frac{m_i}m\right]
$$
and (see the formula for $d(m)$ in Section 3)
$$
  d(m)=n+2\sum_{i=1}^n\left[\frac{m_i}m\right]
$$
(cf. \cite{ser}.)

Making use of the above formula, it is easy to compute the numbers $d(m)$ for all exceptional algebras. They are given in the tables on Fig. 2-6.

\subsection{} A classification of nilpotent elements $e$ in the exceptional Lie algebras follows from the classification of $\sl_2$-triples obtained in \cite{dyn}. The corresponding Levi subalgebras $\l(e)$ and, hence, the Coxeter numbers $h(e)$ also can be derived from the tables of that paper. The centralizers $\z(e)$ were determined in \cite{ela}. The inclusion relation for the closures of nilpotent orbits was described in \cite{spa} (see also \cite{mcg}).

Having all this information, it is easy to find the nilpotent orbits that do not lie in the closure of another nilpotent orbit with the same Coxeter number. The Hasse diagrams for the sets of such orbits are depicted on Fig 2-6, where each involved orbit Ad$(G)e$ is given by the type of $\l(e)'$, and the dimension of $\z(e)$ is indicated in parentheses. For each $m$, the set of orbits with Coxeter number $m$ is situated in the corresponding stripe between dotted lines.

\

\

\

\begin{center}
{\Large\bf G$_2$}

\

\

\begin{tabular}{c|cccc}
$m$&2&3&4&5\\
\hline
$d(m)$&6&4&4&4
\end{tabular}

\

\

\

\begin{tikzpicture}
\node at (-4,1){$m$ ($d(m)$)};%
\node (g2) at  (0,0) [draw,label=right:(2)]  {G$_2$};%
\node at (-4,0) {6 (2)};%
\node (ta1) at   (1.15,-1) [draw,label=right:(6)] {$\widetilde{\textrm A_1}$};%
\node (a1) at   (-1.15,-1) [label=right:(8)] {A$_1$};%
\node at (-4,-1) {2 (6)};%
\node (0) at (0,-2) [draw,label=right:(14)] {0};%
\node at (-4,-2){1 (14)};%
\path
(-5,-.5) edge[dotted] (3,-.5)
(-5,-1.5) edge[dotted] (3,-1.5)
(g2) edge[-stealth] (a1) edge[-stealth] (ta1)
(ta1) edge[-stealth] (0)
(a1) edge[-stealth] (0);
\end{tikzpicture}

\

\

Fig. 2

\pagebreak

\

\vfill

{\Large\bf F$_4$}

\vfill

\begin{tabular}{c|ccccccccccc}
$m$&2&3&4&5&6&7&8&9&10&11&12\\
\hline
$d(m)$&24&16&12&12&8&8&6&6&6&6&4
\end{tabular}

\

\

\begin{tikzpicture}
\node at (-4,1){$m$ ($d(m)$)};%
\node (f4) at  (0,0) [draw,label=right:(4)]  {F$_4$};%
\node at (-4,0) {12 (4)};%
\node (b3) at   (-1.15,-1) [label=right:(10)] {B$_3$};%
\node at (-4,-1){6 (8)};%
\node (c3) at (1.15,-1) [label=right:(10)] {C$_3$};%
\node (b2) at (1.15,-2) [label=right:(16)] {B$_2$};%
\node at (-4,-2){4 (12)};%
\node (a2a1) at (-1.15,-3) [draw,label=right:(16)] {$\widetilde{\textrm A}_2$+A$_1$};%
\node at (-4,-3) {3 (16)};%
\node (a1a1) at (0,-4) [draw,label=right:(24)] {A$_1$+$\widetilde{\textrm A_1}$};%
\node at (-4,-4){2 (24)};%
\node (0) at (0,-5) [draw,label=right:(52)] {0};%
\node at (-4,-5){1 (52)};%
\path
(-5,-.5) edge[dotted] (3,-.5)
(-5,-1.5) edge[dotted] (3,-1.5)
(-5,-2.5) edge[dotted] (3,-2.5)
(-5,-3.5) edge[dotted] (3,-3.5)
(-5,-4.5) edge[dotted] (3,-4.5)
(f4) edge[-stealth] (b3) edge[-stealth] (c3)
(b3) edge[-stealth] (a2a1) edge[-stealth] (b2)
(c3) edge[-stealth] (a2a1) edge[-stealth] (b2)
(a2a1) edge[-stealth] (a1a1)
(b2) edge[-stealth] (a1a1)
(a1a1) edge[-stealth] (0);
\end{tikzpicture}

\vfill

Fig. 3

\

\

\

{\Large\bf E$_6$}

\vfill

\begin{tabular}{c|ccccccccccc}
$m$&2&3&4&5&6&7&8&9&10&11&12\\
\hline
$d(m)$&38&24&20&16&12&12&10&8&8&8&6
\end{tabular}

\

\

\begin{tikzpicture}
\node at (-4,1){$m$ ($d(m)$)};%
\node (e6) at  (0,0) [draw,label=right:(6)] {E$_6$};%
\node at (-4,0){12 (6)};%
\node (d5) at   (0,-1) [draw,label=right:(10)] {D$_5$};%
\node at (-4,-1){8 (10)};%
\node (d4) at (1.5,-2) [label=right:(18)] {D$_4$};%
\node at (-4,-2){6 (12)};%
\node (a4a1) at (-1.5,-3) [draw,label=right:(16)] {A$_4$+A$_1$};%
\node at (-4,-3){5 (16)};%
\node (a3a1) at (0,-4) [label=right:(22)] {A$_3$+A$_1$};%
\node at (-4,-4){4 (20)};%
\node (2a2a1) at (0,-5) [draw,label=right:(24)] {2A$_2$+A$_1$};%
\node at (-4,-5){3 (24)};%
\node (3a1) at (0,-6) [draw,label=right:(38)] {3A$_1$};%
\node at (-4,-6){2 (38)};%
\node (0) at (0,-7) [draw,label=right:(78)] {0};%
\node at (-4,-7){1 (78)};%
\path
(e6) edge[-stealth] (d5)
(d5) edge[-stealth] (a4a1) edge[-stealth] (d4)
(d4) edge[-stealth] (a3a1)
(a4a1) edge[-stealth] (a3a1)
(a3a1) edge[-stealth] (2a2a1)
(2a2a1) edge[-stealth] (3a1)
(3a1) edge[-stealth] (0)
(-5,-.5) edge[dotted] (3,-.5)
(-5,-1.5) edge[dotted] (3,-1.5)
(-5,-2.5) edge[dotted] (3,-2.5)
(-5,-3.5) edge[dotted] (3,-3.5)
(-5,-4.5) edge[dotted] (3,-4.5)
(-5,-5.5) edge[dotted] (3,-5.5)
(-5,-6.5) edge[dotted] (3,-6.5)
;
\end{tikzpicture}

\vfill

Fig. 4

\vfill

\

\pagebreak

\

\vfill

{\Large\bf E$_7$}

\

\

\begin{tabular}{c|cccccccccccccccccc}
$m$&2&3&4&5&6&7&8&9&10&11&12&13&14&15&16&17&18\\
\hline
$d(m)$&63&43&33&27&21&19&17&15&13&13&11&11&9&9&9&9&7
\end{tabular}

\

\

\

\begin{tikzpicture}
\node at (-4,19){$m$ ($d(m)$)};%
\node (e7) at (0,18) [draw,label=right:(7)] {E$_7$};%
\node at (-4,18) {18 (7)};%
\node (e6) at   (-2,17) [label=right:(13)] {E$_6$};%
\node at (-4,17) {12 (11)};%
\node (d6) at (2,16) [label=right:(15)] {D$_6$};%
\node at (-4,16){10 (13)};%
\node (d5a1) at (2,15) [label=right:(19)] {D$_5$+$A_1$};%
\node at (-4,15) {8 (17)};%
\node (a6) at (-2,14) [draw,label=right:(19)] {A$_6$};%
\node at (-4,14){7 (19)};%
\node (a5') at (-2,13) [label=right:(25)] {A$_5'$};%
\node at (-4,13){6 (21)};%
\node (a5a1'') at (2,13) [label=right:(25)] {[A$_5$+A$_1$]$''$};%
\node (d4a1) at (5,13) [label=right:(31)] {D$_4$+A$_1$};%
\node (a4a2) at (0,12) [draw,label=right:(27)] {A$_4$+A$_2$};%
\node at (-4,12){5 (27)};%
\node (a3a2a1) at (0,11) [draw,label=right:(33)] {A$_3$+A$_2$+A$_1$};%
\node at (-4,11){4 (33)};%
\node (2a2a1) at (0,10) [draw,label=right:(43)] {2A$_2$+A$_1$};%
\node at (-4,10){3 (43)};%
\node (4a1) at (0,9) [draw,label=right:(63)] {4A$_1$};%
\node at (-4,9) {2 (63)};%
\node (0) at (0,8) [draw,label=right:(133)] {0};%
\node at (-4,8) {1 (133)};%
\path
(e7) edge[-stealth] (e6) edge[-stealth] (d6)
(e6) edge[-stealth] (a6) edge[-stealth] (d5a1)
(d6) edge[-stealth] (a6) edge[-stealth] (d5a1)
(a6) edge[-stealth] (a5') edge[-stealth] (a5a1'')
(d5a1) edge[-stealth] (a5') edge[-stealth] (a5a1'') edge[-stealth] (d4a1)
(a5') edge[-stealth] (a4a2)
(a5a1'') edge[-stealth] (a4a2)
(a4a2) edge[-stealth] (a3a2a1)
(d4a1) edge[-stealth] (a3a2a1)
(a3a2a1) edge[-stealth] (2a2a1)
(2a2a1) edge[-stealth] (4a1)
(4a1) edge[-stealth] (0)
(-5,17.5) edge[dotted] (6,17.5)
(-5,16.5) edge[dotted] (6,16.5)
(-5,15.5) edge[dotted] (6,15.5)
(-5,14.5) edge[dotted] (6,14.5)
(-5,13.5) edge[dotted] (6,13.5)
(-5,12.5) edge[dotted] (6,12.5)
(-5,11.5) edge[dotted] (6,11.5)
(-5,10.5) edge[dotted] (6,10.5)
(-5,9.5) edge[dotted] (6,9.5)
(-5,8.5) edge[dotted] (6,8.5);
\end{tikzpicture}

\

\

Fig. 5

\vfill

\

\pagebreak

{\Large\bf E$_8$}
\end{center}

\

\

\hskip5em\begin{tabular}{c|cccccccccccccc}
$m$&2&3&4&5&6&7&8&9&10&11&12&13&14&15\\
\hline
$d(m)$&120&80&60&48&40&36&30&28&24&24&20&20&18&16
\end{tabular}

\

\hskip5em\begin{tabular}{c|cccccccccccccccc}
$m$&16&17&18&19&20&21&22&23&24&25&26&27&28&29&30\\
\hline
$d(m)$&16&16&14&14&12&12&12&12&10&10&10&10&10&10&8
\end{tabular}

\

\

\

\begin{center}

\begin{tikzpicture}
\node at (-4,1){$m$ ($d(m)$)};%
\node (e8) at   (0,0) [draw,label=right:(8)]  {E$_8$};%
\node at (-4,0){30 (8)};%
\node (e7) at   (0,-1) [label=right:(16)] {E$_7$};%
\node at (-4,-1){18 (14)};%
\node (d7) at   (1,-2) [label=right:(22)] {D$_7$};%
\node at (-4,-2){12 (20)};%
\node (e6a1) at   (-1,-2) [label=right:(26)] {E$_6$+A$_1$};%
\node (d6a1) at   (1,-3) [label=right:(28)] {D$_6$+A$_1$};%
\node at (-4,-3){10 (24)};%
\node (a7) at (-1,-4) [draw,label=right:(30)] {A$_7$};%
\node at (-4,-4) {8 (30)};%
\node (a6a1) at (0,-5) [draw,label=right:(36)] {A$_6$+A$_1$};%
\node at (-4,-5){7 (36)};%
\node (a5a2) at (0,-6) [label=right:(42)] {A$_5$+A$_2$};%
\node at (-4,-6){6 (40)};%
\node (a4a3) at (0,-7) [draw,label=right:(48)] {A$_4$+A$_3$};%
\node at (-4,-7) {5 (48)};%
\node (2a3) at (0,-8) [draw,label=right:(60)] {A$_3$+A$_3$};%
\node at (-4,-8) {4 (60)};%
\node (2a22a1) at (0,-9) [draw,label=right:(80)] {2A$_2$+2A$_1$};%
\node at (-4,-9) {3 (80)};%
\node (4a1) at (0,-10) [draw,label=right:(120)] {4A$_1$};%
\node at (-4,-10) {2 (120)};%
\node (0) at (0,-11) [draw,label=right:(248)] {0};%
\node at (-4,-11) {1 (248)};%
\path
(-5,-.5) edge[dashed] (3,-.5)
(-5,-1.5) edge[dashed] (3,-1.5)
(-5,-2.5) edge[dashed] (3,-2.5)
(-5,-3.5) edge[dashed] (3,-3.5)
(-5,-4.5) edge[dashed] (3,-4.5)
(-5,-5.5) edge[dashed] (3,-5.5)
(-5,-6.5) edge[dashed] (3,-6.5)
(-5,-7.5) edge[dashed] (3,-7.5)
(-5,-8.5) edge[dashed] (3,-8.5)
(-5,-9.5) edge[dashed] (3,-9.5)
(-5,-10.5) edge[dashed] (3,-10.5)
(e8) edge[-stealth] (e7)
(e7) edge[-stealth] (e6a1) edge[-stealth] (d7)
(d7) edge[-stealth] (d6a1)
(e6a1) edge[-stealth] (a7)
(d6a1) edge[-stealth] (a6a1)
(a7)  edge[-stealth] (a6a1)
(a6a1) edge[-stealth] (a5a2)
(a5a2) edge[-stealth] (a4a3)
(a4a3) edge[-stealth] (2a3)
(2a3) edge[-stealth] (2a22a1)
(2a22a1) edge[-stealth] (4a1)
(4a1) edge[-stealth] (0);
\end{tikzpicture}

\

\

Fig. 6

\end{center}

\

\

\

Theorems 2 and 3 for the exceptional Lie algebras are immediately obtained by observing these diagrams. The nilpotent orbits that turn out to be exceptional are framed there.

\end{document}